\documentclass[twoside,12pt,reqno]{amsart}
\usepackage{amsmath, amsthm,amssymb,amstext,amsfonts,amscd}
\usepackage{graphicx}
\usepackage{caption}
\usepackage{refstyle}
\usepackage{multirow}
\usepackage{hyperref}
\usepackage{lmodern}
\usepackage{latexsym}
\usepackage{amssymb}
\usepackage[utf8]{inputenc}
\usepackage[english]{babel}
\usepackage{wasysym} %\usepackage{fontspec}
\usepackage{cite}
\numberwithin{equation}{section}
\begin{document}

\begin{center}
\textbf{ANALYTICAL EXPRESSION FOR THE EXACT CURVED SURFACE AREA AND VOLUME OF A HYPERBOLOID OF ONE SHEET  VIA MELLIN-BARNES TYPE CONTOUR INTEGRATION}
\end{center}
\begin{center}
\noindent {M.A. Pathan$^{1,2}$, M. I. Qureshi$^3$, Javid Majid$^{3,*}$    }
\end{center}
\begin{center}
$^1$ Centre for Mathematical and Statistical Sciences (CMSS), Peechi,\\
Thrissur-680653, Kerala, India\\
$^2$ Department of Mathematics, Aligarh Muslim University,\\
Aligarh-202002, U.P., India\\
$^3$ Department of Applied Sciences and Humanities \\
Faculty of Engineering and Technology\\
Jamia Millia Islamia (A Central University), New Delhi-110025, India.\\
 Emails: mapathan@gmail.com,~miqureshi\_delhi@yahoo.co.in\\
  $^*$Corresponding author: javidmajid375@gmail.com 
\end{center}
\vskip 0.2 cm
\textbf{Abstract:}
In this article, we aim at obtaining the analytical expression ({\bf not previously found and recorded in the literature}) for the exact curved surface area of a  hyperboloid of one sheet in terms of Srivastava-Daoust triple hypergeometric function. The derivation is based on Mellin-Barnes type contour integral representations of generalized hypergeometric function$~_pF_q(z)$, Meijer's $G$-function, decomposition formula for Meijer's $G$-function and series rearrangement technique. Further, we also obtain the formula for the volume of a hyperboloid of one sheet. The closed forms for the exact curved surface area and volume of the hyperboloid of one sheet are also verified numerically by using {\it Mathematica Program}.
\vskip 0.2 cm
\section{\bf{Introduction and preliminaries}}
\noindent
For the definition of Pochhammer symbols, power series form of generalized hypergeometric function ${}_pF_q(z)$ and several related results, we refer the beautiful monographs (see, e.g., \cite{And1999, Erd1, Lebedev1972,Luke169 , Rain, Slater1966, SriMano})

\noindent
$\bullet$ Gauss' classical summation theorem\cite[p.49, Th.(18)]{Rain} is given by:\\
\indent
\begin{equation}\label{A2.2}
~\mathbf~{_{2}F_{1}}\left[\begin{array}{ccc}
                                        \alpha,~~\beta;   \\
                                        & 1\\
                                        ~~~~~~\gamma; \\
                                        \end{array}\right]
=\frac{\Gamma (\gamma)~\Gamma (\gamma-\alpha-\beta)}{\Gamma(\gamma-\alpha)\Gamma(\gamma-\beta)},~~~~~~~~~~~~~~~~~~~~~~~~~~~~
\end{equation}
where $\mathfrak{Re}(\gamma-\alpha-\beta)>0$ and $\gamma\in \mathbb{C}\backslash\mathbb Z_{0}^-$.\\

\noindent
$\bullet$ Analytic continuation formula \cite[p.559, Entry(15.3.9),\cite{Erd1}, p.107, Eq.(33), p.105, Eq.(7) and p.106, Eq.(23), \cite{Lebedev1972}, p.250, Eq.(9.5.10), \cite{MagOberSoni1966}, p.48, Eq.(3rd),  \cite{Prudbrimari90}, p.454, Entry(8)]{Abramowitz}:
\begin{equation*}
~_2F_1\left[\begin{array}{ll}
A,B;\\
& z\\
~~~C;\end{array}\right]=\frac{\Gamma(C)\Gamma(C-A-B)}{\Gamma(C-A)\Gamma(C-B)}z^{-A}~_2F_1\left[\begin{array}{ll}
~~~A,1+A-C;\\
& 1-\frac{1}{z}\\
A+B-C+1;\end{array}\right]+
\end{equation*}
\begin{equation}\label{A15.30}
+\frac{\Gamma(C)\Gamma(A+B-C)}{\Gamma(A)\Gamma(B)}z^{A-C}(1-z)^{C-A-B}~_2F_1\left[\begin{array}{ll}
~~~C-A,1-A;\\
& 1-\frac{1}{z}\\
1+C-A-B;\end{array}\right]
\end{equation}
where $|\arg(1-z)|< \pi,|\arg(z)|< \pi,|z|>1,\mathfrak{R}(z)>\frac{1}{2} $ and $A+B-C\neq 0,\pm1,\pm2,....$

\noindent
$\bullet$ Mellin-Barnes type contour integral representation of binomial function$~{}_1F_0(z)$:
\begin{equation}\label{D12.100}
(1-z)^{-a}=~{}_1F_0\left[\begin{array}{ll}
~a;\\
&z\\
-;\end{array}\right]=\frac{1}{2\pi i~\Gamma(a)}\int_{-i\infty}^{+i\infty}\Gamma(a+s)\Gamma(-s)(-z)^s~ds:~~z\neq 0,
\end{equation}
where $|\arg(-z)|<\pi, %|z|<1,
a\in\mathbb{C}\backslash\mathbb{Z}^-_0$ and $i=\sqrt{(-1)}$.\\

\noindent
$\bullet$ Mellin-Barnes type contour integral representation of $~{}_pF_q(z)$ \cite[p.43, Eq.(6)]{SriMano}:
\begin{equation}\label{D12.101}
~{}_pF_q\left[\begin{array}{ll}
\alpha_1,\alpha_2,...,\alpha_p;\\
& z\\
\beta_1,\beta_2,...,\beta_q;\end{array}\right]=\frac{\Gamma(\beta_1)\Gamma(\beta_2)...\Gamma(\beta_q)}{\Gamma(\alpha_1)\Gamma(\alpha_2)...\Gamma(\alpha_p)}.\frac{1}{2\pi i}\int_{\sigma-i\infty}^{\sigma+i\infty}\frac{\Gamma(\alpha_1+\xi)...\Gamma(\alpha_p+\xi)\Gamma(-\xi)(-z)^\xi}{\Gamma(\beta_1+\xi)...\Gamma(\beta_q+\xi)}~d\xi,
\end{equation}
where $z\neq0$.\\

\noindent
 Convergence conditions:\\
If $p=q+1$, then $|\arg(-z)|<\pi$.\\
 %or $|z|<1$,\\
If $p=q$, then $|\arg(-z)|<\frac{\pi}{2}$,\\
% or $\Re e(z)<0$ \\
and no $\alpha_i~(i=1,2,...,p)$ %nor $\beta_j$
is zero or negative integer but some of $\beta_j~(j=1,2,...,q)$ may be zero or negative integers.\\

\noindent
$\bullet$ Mellin-Barnes type contour integral representation of Meijer's G-function \cite[p.45, Eq.(1)]{SriMano}, see also\cite{Erd1,Luke169}:\\
When $p\leq q$ and $1\leq m \leq q,~0\leq n\leq p,$ then
\begin{equation*}
G^{m,n}_{p,q}\left(z\left|   \begin{array}{ll}
\alpha_1,\alpha_2,\alpha_3,...,\alpha_n;\alpha_{n+1},...,\alpha_p\\
\beta_1,\beta_2,\beta_3,...,\beta_m;\beta_{m+1},...,\beta_q  \end{array}\right. \right)=~~~~~~~~~~~~~~~~~~~~~~~~~~~~~~
\end{equation*}
\begin{equation*}
=\frac{1}{2\pi i}\int_{-i\infty}^{+i\infty}\frac{\prod_{j=1}^{m}\Gamma(\beta_j-s)\prod_{j=1}^{n}\Gamma(1-\alpha_j+s)}{\prod_{j=m+1}^{q}\Gamma(1-\beta_j+s)\prod_{j=n+1}^{p}\Gamma(\alpha_j-s)}(z)^s~ds
~~~~~~~~
\end{equation*}
\begin{equation}\label{E12.12}
=\frac{1}{2\pi i}\int_{-i\infty}^{+i\infty}\frac{\Gamma(\beta_1-s)...\Gamma(\beta_m-s)\Gamma(1-\alpha_1+s)...\Gamma(1-\alpha_n+s)}{\Gamma(1-\beta_{m+1}+s)...\Gamma(1-\beta_q+s)\Gamma(\alpha_{n+1}-s)...\Gamma(\alpha_p-s)}(z)^s~ds,
\end{equation}
where $z\neq 0,(\alpha_i-\beta_j)\neq$ positive integers, $i=1,2,3,...,n;~j=1,2,3,...,m$. For details of contours, see\cite[p.207, \cite{Luke169}, p.144]{Erd1}.\\
\vskip 0.2cm
\noindent
Convergence conditions of Meijer's G-function:\\
When $\Lambda=m+n-\left(\frac{p+q}{2}\right),~\nu=\sum_{j=1}^{q}\beta_j-\sum_{j=1}^{p}\alpha_j, $ then
\begin{enumerate}

\item[{(i)}] The integral (\ref{E12.12}) is convergent when $|arg(z)|<\Lambda\pi$ and $\Lambda> 0$.\\

\item[{(ii)}] If $|arg(z)|=\Lambda\pi$ and $~\Lambda\geq 0$, then the integral (\ref{E12.12}) is absolutely convergent when $p=q$ and $\mathfrak{R}(\nu)<-1$.

\item[{(iii)}]If $|arg(z)|=\Lambda\pi$ and $~\Lambda\geq 0$, then the integral (\ref{E12.12}) is also absolutely convergent, when $p\neq q,~(q-p)\sigma>\mathfrak{R}(\nu)+1-\left(\frac{q-p}{2}\right)$ and $s=\sigma+ik$, where $\sigma$ and $k$ are real. $\sigma$ is chosen so that for $k  \rightarrow\pm\infty$.\\

For other two types of contours, following will be convergence conditions of the integral (\ref{E12.12}) \\

\item[{(iv)}]The integral (\ref{E12.12}) is  convergent if $q\geq 1$ and either $p<q,0<|z|<\infty$ or $p=q,~0<|z|<1$.

\item[{(v)}]The integral (\ref{E12.12}) is  convergent if $p \geq 1$ and either $p>q,0<|z|<\infty$ or $p=q,~|z|>1$.
\end{enumerate}

\noindent
$\bullet$ Relations between Meijer's $G$- function and ${}_2F_1(z)$ \cite[p.61, \cite{Wikipedia}, p.77, Eq.(1)]{Mathaisaxena}:
\begin{equation*}
G^{2~2}_{2~2}\left(z\left|   \begin{array}{ll}
1-a,1-b;-\\
0,c-a-b;-  \end{array}\right. \right)=\frac{\Gamma(a)\Gamma(b)\Gamma(c-a)\Gamma(c-b)}{\Gamma(c)}{}_2F_1\left[\begin{array}{ll}
a,~b;\\
& 1-z\\
~~~ ~c;\end{array}\right],
\end{equation*}

where $|1-z|<1$ and $c-a,c-b\neq 0,-1,-2,...$
\begin{equation*}
G^{2~2}_{2~2}\left(z\left|   \begin{array}{ll}
a_1,a_2;-\\
b_1,b_2;-  \end{array}\right. \right)=\frac{\Gamma(1-a_1+b_1)\Gamma(1-a_1+b_2)\Gamma(1-a_2+b_1)\Gamma(1-a_2+b_2)z^{b_1}}{\Gamma(2-a_1-a_2+b_1+b_2)}\times
\end{equation*}
\begin{equation}\label{F12.200}
\times{}_2F_1\left[\begin{array}{ll}
1-a_1+b_1,1-a_2+b_1;\\
& 1-z\\
~~~2-a_1-a_2+b_1+b_2;\end{array}\right];~~~|1-z|<1.\\
\end{equation}
\noindent

\noindent
$\bullet$ Decomposition formula for Meijer's $G$- function \cite[p.208, \cite{Wikipedia}]{Erd1}:\\

\noindent
If no two $a_k,~k=1,2,...,n$, differ by an integer or zero and $p\geq 1$, then
\begin{equation*}
G^{m,n}_{p,q}\left(z\left|   \begin{array}{ll}
a_1,...,a_n;a_{n+1},...,a_p\\
b_1,...,b_m;b_{m+1},...,b_q  \end{array}\right. \right)==\sum_{h=1}^{n}\frac{\prod_{\stackrel{j=1}{j\neq h}}^n\Gamma(
a_h-a_j)\prod_{j=1}^{m}\Gamma(1+b_j-a_h)}{\prod_{j=n+1}^{p}\Gamma(1+a_j-a_h)\prod_{j=m+1}^{q}\Gamma(a_h-b_j)}z^{({a_h}-1)}\times
\end{equation*}

\begin{equation}\label{A15.21}
\times{}_qF_{p-1}\left[\begin{array}{ll}
~~~~~~~~~~~~~~~~~~~~~~~~~~~~~~~~~~~~~~1+b_1-a_h~,...,~1+b_q-a_h~~;\\
& \frac{(-1)^{(q-m-n)}}{z}\\
1+a_1-a_h,...,1+a_{h-1}-a_h,1+a_{h+1}-a_h,...,1+a_p-a_h;\end{array}\right],
\end{equation}

Convergence conditions:\\
When $(a_j-a_h)\notin\mathbb{Z};1\leq j\leq n;1\leq h\leq n$, then
\begin{enumerate}

\item[{(i)}] The formula (\ref{A15.21}) is convergent when  $p> q$.
\item[{(ii)}]The formula (\ref{A15.21}) is convergent when  $p=q$ and $m+n=p+1$ and $z\notin(-1,0)$.

\item[{(iii)}]The formula (\ref{A15.21}) is convergent when  $p=q$ and $m+n>p+1$.

\item[{(iv)}]The formula (\ref{A15.21}) is convergent when  $p=q$ and $m+n=p$ and $|z| >1$.\\
\end{enumerate}

%where $q<p$ or $p=q$ and $|z|>1$. Here the prime in $\prod^{'}$ denotes the omission of the factor $\Gamma(a_h-a_h)$ and the asterisk in$~_pF_q$ indicates the omission of the parameter $(1+a_h-a_h)$.\\

\noindent
$\bullet$ The following generalization of the hypergeometric function in several variables has been given by Srivastava and Daoust ( \cite[pp.199-200, Eq.(2.1),\cite{Sridaoust1972}]{Sridaoust1969202}) which is referred to, in the literature as the generalized Lauricella function of several variables (see also \cite[ p.454, Eq.(4.1)]{Sridaoust1969457}):

\begin{equation*}
 F_{C:~D^{(1)};...;~D^{(n)}}^{A:~B^{(1)};...;~B^{(n)}}
\begin{pmatrix}
[(a_A):~\vartheta^{(1)},...,~\vartheta^{(n)}]:[(b^{(1)}_{B^{(1)}}):\varphi^{(1)}];...;[(b^{(n)}_{B^{(n)}}):\varphi^{(n)}];\\
&   x_1,...,x_n \\
[(c_C)~:~\psi^{(1)},...,~\psi^{(n)}]:[(d^{(1)}_{D^{(1)}}):\delta^{(1)}];...;[(d^{(n)}_{D^{(n)}}):\delta^{(n)}];\\
\end{pmatrix}
\end{equation*}

\begin{equation}\label{AB15.42}
=\sum_{m_1,m_2,...,m_n=0}^{\infty}\Omega(m_1,m_2,...,m_n)~\frac{x_1^{m_1}}{m_1!}\frac{x_2^{m_2}}{m_2!}\cdots\frac{x_n^{m_n}}{m_n!},
\end{equation}
where
\begin{equation}
\Omega(m_1,m_2,...,m_n):=\frac{\prod_{j=1}^A(a_j)_{{m_1}\vartheta^{(1)}_j+\cdots+{m_n}\vartheta^{(n)}_j}\prod_{j=1}^{B^{(1)}}(b^{(1)}_j)_{{m_1}\varphi^{(1)}_j}...~\prod_{j=1}^{B^{(n)}}(b^{(n)}_j)_{{m_n}\varphi^{(n)}_j}}{\prod_{j=1}^C(c_j)_{{m_1}\psi^{(1)}_j+\cdots+{m_n}\psi^{(n)}_j}\prod_{j=1}^{D^{(1)}}(d^{(1)}_j)_{{m_1}\delta^{(1)}_j}...~\prod_{j=1}^{D^{(n)}}(d^{(n)}_j)_{{m_n}\delta^{(n)}_j}}
\end{equation}

and the coefficients
$$\vartheta_j^{(i)}, j=1,2,..., A;~~\varphi_j^{(i)}, j=1,2,..., B^{(i)};~~\psi_j^{(i)}, j=1,2,..., C;~~\delta_j^{(i)}, j=1,2,..., D^{(i)};$$
for all $i \in \big\{1,2,..., n\big\}$ are real and positive,\\
then, with the positive constants $\vartheta's,\varphi's,\psi's$ and $\delta's$ equated to one,\\

$F_{0:1;...;1}^{1:1;...;1}$~
will correspond to Lauricella's $F_A^{(n)}$-function,
\vskip 0.2cm
$F_{1:0;...;0}^{0:2;...;2}~$
will correspond to Lauricella's $F_B^{(n)}$-function,
\vskip 0.2cm
$F_{0:1;...;1}^{2:0;...;0}$~ will correspond to Lauricella's $F_C^{(n)}$-function and
\vskip 0.1cm
$F_{1:0;...;0}^{1:1;...;1}$~ will correspond to Lauricella's fourth function $F_D^{(n)}$  \cite[p.113]{Lauricella93158}; \\
\vskip 0.2cm
while
$F_{0:D^{(1)};...;D^{(n)}}^{0:B^{(1)};...;B^{(n)}}$~
will yield the product
\begin{equation}
_{B^{(1)}}{F}_{D^{(1)}}\left[\begin{array}{ccc}
(b_{B^{(1)}}^{(1)});\\
& x_1\\
(d_{D^{(1)}}^{(1)});\\
\end{array}\right]..._{B^{(n)}}{F}_{D^{(n)}}\left[\begin{array}{ccc}
(b_{B^{(n)}}^{(n)});\\
& x_n\\
(d_{D^{(n)}}^{(n)});\\
\end{array}\right]
\end{equation}
of $n$ generalized hypergeometric functions with different arguments.\\
The multiple hypergeometric function (\ref{AB15.42}) is the generalization of Fox-Wright hypergeometric function of one variable $_p{\Psi}_q$ and $_p{\Psi}_q^*$ \cite{Wright1935,Wright1940}.

Let
\begin{equation}
E_i=\left( \mu_i^{1+\sum_{j=1}^{D^{(i)}}\delta_j^{(i)}-\sum_{j=1}^{B^{(i)}}\varphi_j^{(i)}}\right) ~\frac{\prod_{j=1}^{C}\left(\sum_{\ell=1}^{n}\mu_{\ell}~\psi_{j}^{(\ell)} \right)^{\psi_{j}^{(i)}}~\prod_{j=1}^{D^{(i)}}(\delta_j^{(i)})^{\delta_j^{(i)}}}{\prod_{j=1}^{A}\left(\sum_{\ell=1}^{n}\mu_{\ell}~ \vartheta_{j}^{(\ell)} \right)^{\vartheta_{j}^{(i)}}~\prod_{j=1}^{B^{(i)}}(\varphi_j^{(i)})^{\varphi_j^{(i)}}},
\end{equation}

\begin{equation}
\varDelta_i= 1+\sum_{j=1}^{C}\psi_j^{(i)} + \sum_{j=1}^{D^{(i)}} \delta_j^{(i)} -\sum_{j=1}^{A}\vartheta_j^{(i)} - \sum_{j=1}^{B^{(i)}}\varphi_j^{(i)};~i=1,2,...,n.~~~~~~~~~~~~~
\end{equation}
\textbf{Case I.} The multiple power series in (\ref{AB15.42}) is convergent for all finite complex values or real values of $x_1,x_2,...,x_n$, when $\varDelta_i>0$, $i=1,2,...,n.$\\
\textbf{Case II.} The multiple power series in (\ref{AB15.42}) is convergent when $\varDelta_1=\varDelta_2=...=\varDelta_n=0$; $|x_1|<\varrho_1,|x_2|<\varrho_2,...,|x_n|<\varrho_n;\\$

where
\begin{equation}
\varrho_i=\min_{\mu_1 ,...,\mu_n>0}\{E_i\},~~~i=1,2,...,n.~~~~~~~~~~~~~~~~~~~~~~~~~~~~~~~~~~
\end{equation}
\textbf{Case III.} The multiple power series in (\ref{AB15.42}) would diverge except when, trivially, $x_1=x_2=...=x_n=0$ when  $\varDelta_i<0$, $i=1,2,...,n.$\\
\begin{center}
	\textbf{Further analysis of Case II}
\end{center}

When
\begin{equation}
\vartheta_j^{(1)}=\vartheta_j^{(2)}=...=\vartheta_j^{(n)}=\vartheta_j,~1\leq j\leq A,\\~~~~~~~~~~~~~~~~~~~~
\end{equation}
\begin{equation} \psi_j^{(1)}=\psi_j^{(2)}=...=\psi_j^{(n)}=\psi_j,~1\leq j\leq C,~~~~~~~~~~~~~~~~~~~
\end{equation}
\begin{equation}
G_i=\frac{\prod_{j=1}^{C}(\psi_{j}^{(i)})^{\psi_{j}^{(i)}} ~\prod_{j=1}^{D^{(i)}}(\delta_j^{(i)})^{\delta_j^{(i)}}}{\prod_{j=1}^{A}(\vartheta_{j}^{(i)})^{\vartheta_{j}^{(i)}} ~\prod_{j=1}^{B^{(i)}}(\varphi_j^{(i)})^{\varphi_j^{(i)}}},~~i=1,2,...,n,~~~~~~~~~
\end{equation}
\begin{equation}
\mho_i\equiv 1+\sum_{j=1}^{C}\psi_j^{(i)} + \sum_{j=1}^{D^{(i)}} \delta_j^{(i)} -\sum_{j=1}^{A}\vartheta_j^{(i)} - \sum_{j=1}^{B^{(i)}}\varphi_j^{(i)},~~i=1,2,...,n
\end{equation}

and
\begin{equation}
\Omega=\sum_{j=1}^{A}\vartheta_j - \sum_{j=1}^{C}\psi_j.~~~~~~~~~~~~~~~~~~~~~~~~~~~~~~~~~~~~~~~~~~~~~~~~~~~~~
\end{equation}
\noindent
\textbf{Case II(a).} The multiple power series in (\ref{AB15.42}) is convergent when ~$\mho_1=\mho_2=...=\mho_n=0;~\Omega>0$ and\\
\begin{equation}
\left( \frac{|x_1|}{G_1}\right) ^{\frac{1}{\Omega}}+...+\left( \frac{|x_n|}{G_n}\right) ^{\frac{1}{\Omega}}<1.~~~~~~~~~~~~~~~~~~~~~~~~~~~~~~~~~~~~~~~~
\end{equation}
\textbf{Case II(b).} The multiple power series in (\ref{AB15.42}) is convergent when  $\mho_1=\mho_2=...=\mho_n=0$; $\Omega\leq 0$ and\\
\begin{equation}
\max \left( \frac{|x_1|}{G_1},\frac{|x_2|}{G_2},...,\frac{|x_n|}{G_n}\right)<1.~~~~~~~~~~~~~~~~~~~~~~~~~~~~~~~~~~~~~
\end{equation}

\noindent
$\bullet$Suppose $\phi(x,y)=0$ is the projection of the curved surface of three dimensional figure $z=f(x,y)$ over the $x$-$y$ plane, then curved surface area is given by
\begin{equation}\label{a12.9}
\hat{S}=\underbrace{\int\int}_{\stackrel{\text{over the area }} {\phi(x,y)=0}}\sqrt{\left\lbrace 1+\left( \frac{\partial z}{\partial x}\right) ^2+\left( \frac{\partial z}{\partial y}\right) ^2\right\rbrace }~dx~dy.
\end{equation}
\begin{figure}[ht!]
\noindent\begin{minipage}{0.60\textwidth}% adapt widths of minipages to your needs
\includegraphics[width=\linewidth]{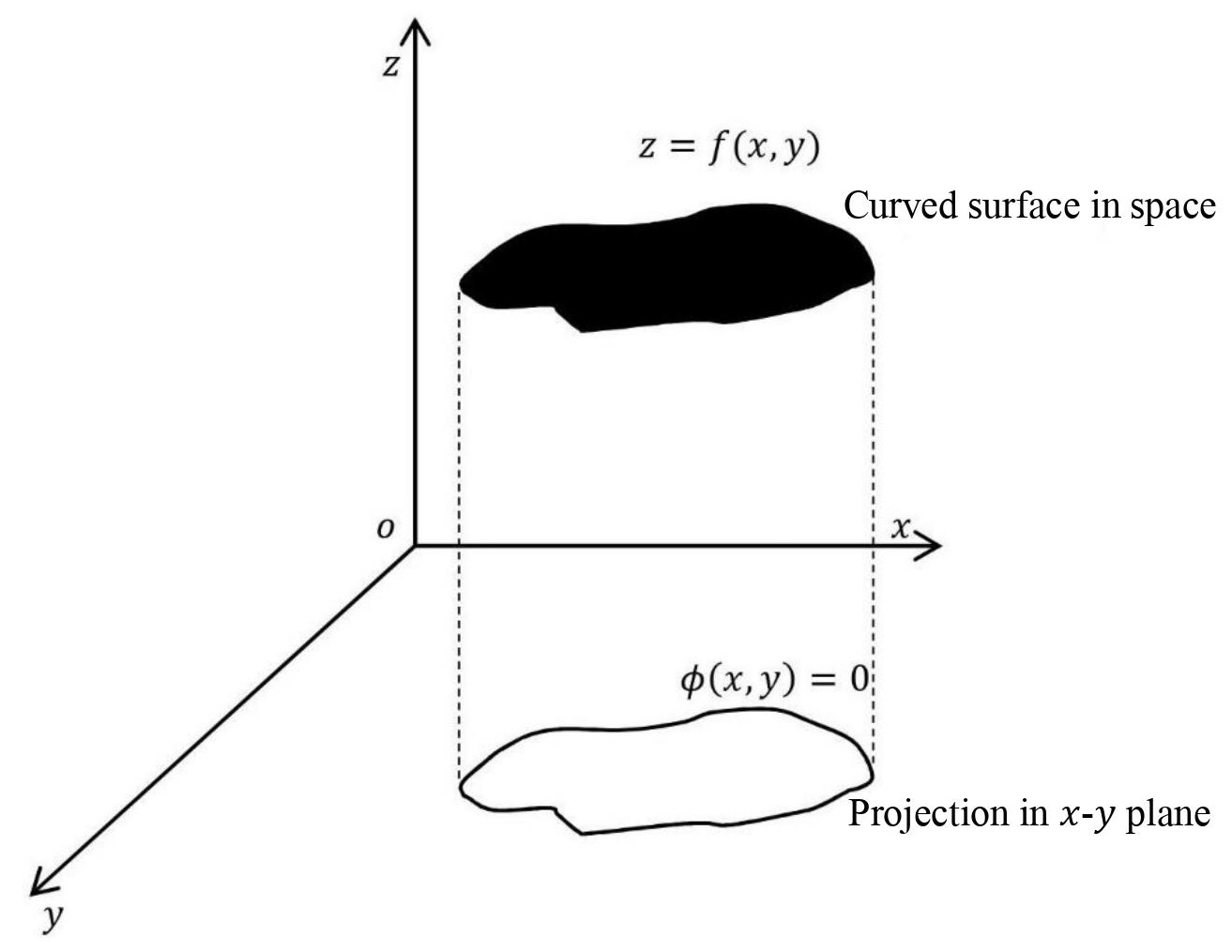}
\end{minipage}%
\caption[short]{Projection of curved surface in $x$-$y$ plane.}
\label{fig:Projection }
\end{figure}

\noindent
\begin{equation}\label{D12.1000}
\bullet\text{A definite integral }\int_{\theta=0}^{\frac{\pi}{2}}\sin^\alpha {\theta}\cos^\beta {\theta}~d\theta=\frac{\Gamma\left(\frac{\alpha+1}{2}\right) \Gamma\left(\frac{\beta+1}{2}\right) }{2\Gamma\left(\frac{\alpha+\beta+2}{2}\right) },~~~~~~~~~~~~~~
\end{equation}

where $\mathfrak{R}(\alpha)>-1,~\mathfrak{R}(\beta)>-1$.\\

The present article is organized as follows. In {\bf section 2}, we evaluated some important definite integrals $\int_{\theta=-\pi}^{\pi}\left(\frac{\cos^2{\theta}}{\sigma^2}+\frac{\sin^2{\theta}}{\lambda^2}\right)^sd\theta$ and $\int_{r=1}^{\lambda} \frac{r^{2s+1}}{(1-r^2)^s}dr$ with suitable convergence conditions using Mellin-Barnes type contour integral representations of binomial function$~_1F_0(z)$, Meijer's G-function, analytic continuation formula for Gauss function and series manipulation technique. These integrals are useful and help us in the derivation of closed form for the exact curved surface area of a hyperboloid of one sheet. In {\bf section 3}, we derive the closed form for obtaining the exact curved surface area of a hyperboloid of one sheet by using Mellin-Barnes type contour integral representations of generalized hypergeometric function$~_pF_q(z)$, Meijer's $G$-function and decomposition formula for Meijer's $G$-function; in terms of generalized hypergeometric function of Srivastava and Daoust. In {\bf section 4}, we derive the formula for the volume of a hyperboloid of one sheet.

\section{\bf Evaluation of some useful definite integrals}
\noindent
The following definite integrals hold true associated with suitable convergence conditions:
\begin{equation}\label{G12.10}
\text{{\bf Theorem 1.}}~\int_{\theta=-\pi}^{\pi}\left(\frac{\cos^2{\theta}}{\sigma^2}+\frac{\sin^2{\theta}}{\lambda^2}\right)^s~d\theta =\frac{2\pi \lambda}{\sigma^{1+2s}}~{}_2F_1\left[\begin{array}{ll}
\frac{1}{2},1+s;\\
&1-\frac{\lambda^2}{\sigma^2}\\
~~~~~~~~1;\end{array}\right],~~~~~~~~~~~~~~~
\end{equation}
where $\sigma\geq \lambda>0$ and it is obvious that $0\leq(1-\frac{\lambda^2}{\sigma^2})<1$.\\
\begin{equation}\label{G12.10000}
\text{{\bf Theorem 2.}}~\int_{\theta=-\pi}^{\pi}\left(\frac{\cos^2{\theta}}{\sigma^2}+\frac{\sin^2{\theta}}{\lambda^2}\right)^s~d\theta =\frac{2\pi \sigma}{\lambda^{1+2s}}~{}_2F_1\left[\begin{array}{ll}
\frac{1}{2},1+s;\\
&1-\frac{\sigma^2}{\lambda^2}\\
~~~~~~~~1;\end{array}\right],~~~~~~~~~~~~~~~
\end{equation}
where $\lambda\geq \sigma>0$ and it is obvious that $0\leq(1-\frac{\sigma^2}{\lambda^2})<1$.\\

\noindent
{\bf Theorem 3.} When $1\leq\lambda<\infty$, then
\begin{equation}\label{A15.22}
~\int_{r=1}^{\lambda} \frac{r^{2s+1}}{(1-r^2)^s}~~dr=\frac{\lambda^{2s}(1-\lambda^2)^{1-s}~\Gamma(s-1)}{2\lambda^2~ \Gamma(s)}~{}_2F_1\left[\begin{array}{ll}
2,1-s;\\
& 1-\frac{1}{\lambda^2}\\
~~~2-s;\end{array}\right],
\end{equation}
where $\mathfrak{R}(s)<1$. \\
\noindent
{\em{Remark:}} The above formulas (\ref{G12.10}), (\ref{G12.10000}) and (\ref{A15.22}) are also verified numerically using Mathematica program.\\

\noindent
{\bf Independent demonstration of the assertions (\ref{G12.10}) and (\ref{G12.10000}) }
\begin{equation*}
\text{Suppose} ~I_1=\int_{\theta=-\pi}^{\pi}\left(\frac{\cos^2{\theta}}{\sigma^2}+\frac{\sin^2{\theta}}{\lambda^2}\right)^s~d\theta.
\end{equation*}
Employing the formulas (\ref{D12.1001}) and (\ref{D12.1002}), we get
\begin{equation*}
I_1=\frac{4}{\lambda^{2s}}\int_{\theta=0}^{\frac{\pi}{2}}\left( \sin^2 \theta\right)^s\left\lbrace 1+\frac{\lambda^2\cos^2{\theta}}{\sigma^2\sin^2{\theta}}\right\rbrace ^s~d\theta
\end{equation*}
\begin{equation}\label{G12.100}
=\frac{4}{\lambda^{2s}}\int_{\theta=0}^{\frac{\pi}{2}} \sin^{2s} \theta~{}_1F_0\left[\begin{array}{ll}
~~-s;\\
&\frac{-\lambda^2\cos^2{\theta}}{\sigma^2\sin^2{\theta}}\\
~-~~;\end{array}\right]~d\theta.
\end{equation}
Employing the contour integral (\ref{D12.100}) of ${}_{1}F_{0}(.)$, we get
\begin{equation}\label{G12.1005}
I_1=\frac{2}{\pi i \Gamma(-s)~\lambda^{2s}}\int_{\theta=0}^{\frac{\pi}{2}} \sin^{2s}\left\lbrace \int_{\zeta=-i\infty}^{+i\infty}\Gamma(-\zeta)\Gamma(-s+\zeta)\left(\frac{\lambda^2\cos^2{\theta}}{\sigma^2\sin^2{\theta}}\right)^\zeta~d\zeta \right\rbrace  ~d\theta.
\end{equation}
Interchanging the order of integration in double integral of (\ref{G12.1005}), we get
\begin{equation}
I_1=\frac{2}{\pi i \Gamma(-s)~\lambda^{2s}}\int_{\zeta=-i\infty}^{+i\infty}\Gamma(-\zeta)\Gamma(-s+\zeta)\left(\frac{\lambda^2}{\sigma^2}\right)^\zeta\left\lbrace \int_{\theta=0}^{\frac{\pi}{2}}\sin^{2s-2\zeta}{\theta}~\cos^{2\zeta}{\theta}~d\theta\right\rbrace ~d\zeta.
\end{equation}
Using the Wallis integral formula (\ref{D12.1000}), we get

\begin{equation*}
I_1=\frac{1}{\pi i \Gamma(-s)\Gamma(1+s)~\lambda^{2s}}\int_{\zeta=-i\infty}^{+i\infty}\Gamma(-\zeta)\Gamma(-s+\zeta)\Gamma\left(\frac{1}{2}+s-\zeta\right)\Gamma\left(\frac{1}{2}+\zeta\right)\left(\frac{\lambda^2}{\sigma^2}\right)^\zeta~d\zeta
\end{equation*}
\begin{equation}
=\frac{1}{\pi i \Gamma(-s)\Gamma(1+s)~\lambda^{2s}}\int_{\zeta=-i\infty}^{+i\infty}\Gamma(0-\zeta)\Gamma(1-(s+1)+\zeta)\Gamma\left(\frac{1}{2}+s-\zeta\right)\Gamma\left(1-\frac{1}{2}+\zeta\right)\left(\frac{\lambda^2}{\sigma^2}\right)^\zeta~d\zeta.
\end{equation}
Applying the definition (\ref{E12.12}) of Meijer's $G$-function, we get
\begin{equation}\label{G12.101}
I_1=\frac{2}{\Gamma(-s)\Gamma(1+s)~\lambda^{2s} }~G^{2~2}_{2~2}\left(\frac{\lambda^2}{\sigma^2}\left|   \begin{array}{ll}
s+1,\frac{1}{2};-\\
0,\frac{1}{2}+s;-  \end{array}\right. \right).~~~~~~~~~~~~~~~~~~~~~~~~~~~~~~~~~~~~
\end{equation}
Employing the conversion formula (\ref{F12.200}) in equation (\ref{G12.101}), and after further simplification, we arrive at the result (\ref{G12.10}).\\

\noindent
The proof of the result (\ref{G12.10000}) follows the same steps as in the proof of (\ref{G12.10}). So we omit the details here.
\vskip 0.2cm
\noindent
{\bf Independent demonstration of the assertion (\ref{A15.22}) }\\

When $1\leq r<\lambda<\infty$, then
\begin{equation*}
\text{Suppose}~I_2=\int_{r=1}^{r=\lambda} \frac{r^{2s+1}}{(1-r^2)^s}~~dr=\int_{1}^{\lambda}r^{2s+1}{(1-r^2)^{-s}}~~dr~~~~~~~~~~~~~~~~~~~~~~~~~~~~~~~~~~~~~~~~~~~`
\end{equation*}

\begin{equation}
=\int_{1}^{\lambda} r^{2s+1}~{}_1F_0\left[\begin{array}{ll}
~s~;\\
& r^2\\
~-;\end{array}\right]~dr.~~~~~~~~~~~~
\end{equation}
Employing the contour integral (\ref{D12.100}) of ${}_{1}F_{0}(.)$, we get

\begin{equation}\label{A15.40}
I_2=\int_{1}^{\lambda} \frac{r^{2s+1}}{2\pi i \Gamma(s)}\left\lbrace \int_{\zeta=-i\infty}^{+i\infty}\Gamma(s+\zeta)\Gamma(-\zeta)\left(-r^2\right)^\zeta~d\zeta \right\rbrace  ~dr.
\end{equation}
Interchanging the order of integration in double integral of (\ref{G12.1005}), we get
\begin{equation*}
I_2=\frac{1}{2\pi i \Gamma(s)}\int_{\zeta=-i\infty}^{+i\infty}\Gamma(s+\zeta)\Gamma(-\zeta)\left(-1\right)^\zeta\left\lbrace \int_{1}^{\lambda}r^{2s+2\zeta+1}~dr\right\rbrace ~d\zeta
\end{equation*}
\begin{equation*}
=\frac{1}{4\pi i \Gamma(s)}\int_{\zeta=-i\infty}^{+i\infty}\Gamma(s+\zeta)\Gamma(-\zeta)\left(-1\right)^\zeta \frac{ (\lambda^{2s+2\zeta+2}-1)}{(1+s+\zeta)}  ~d\zeta
\end{equation*}

\begin{equation*}
=\frac{\lambda^{2s+2}}{4\pi i \Gamma(s)}\int_{\zeta=-i\infty}^{+i\infty}\frac{\Gamma(s+\zeta)\Gamma(1+s+\zeta)\Gamma(-\zeta)\left(-\lambda^{2}\right)^\zeta ~d\zeta}{\Gamma(2+s+\zeta)}-
\end{equation*}
\begin{equation}\label{A15.23}
-\frac{1}{4\pi i \Gamma(s)}\int_{\zeta=-i\infty}^{+i\infty}\frac{\Gamma(s+\zeta)\Gamma(1+s+\zeta)\Gamma(-\zeta)\left(-1\right)^\zeta ~d\zeta}{\Gamma(2+s+\zeta)}.
\end{equation}

Applying contour integral (\ref{D12.101}) of $~_pF_q(z)$ in both contour integrals of right hand side of equation (\ref{A15.23}), we get
\begin{equation}\label{A15.41}
I_2=\frac{\lambda^{2s+2}}{2(1+s)}~_2F_1\left[\begin{array}{ll}
s,1+s;\\
& \lambda^2\\
~~2+s;\end{array}\right]-\frac{1}{2(1+s)}~_2F_1\left[\begin{array}{ll}
s,1+s;\\
& 1\\
~~2+s;\end{array}\right],
\end{equation}
where $|\lambda^2| \geq1$ and $\mathfrak{R}(s)<1$.\\

\noindent
Employing the analytic continuation formula (\ref{A15.30}) and Gauss classical summation theorem (\ref{A2.2}) in equation (\ref{A15.41}), we get
\begin{equation*}
I_2=\frac{\lambda^{2}\Gamma(1+s)\Gamma(1-s)}{2}~_1F_0\left[\begin{array}{ll}
-1;\\
& 1-\frac{1}{\lambda^2}\\
~-;\end{array}\right]+\frac{\lambda^{2s}(1-\lambda^2)^{1-s}~\Gamma(s-1)}{2\lambda^2~\Gamma(s)}\times
\end{equation*}

\begin{equation}\label{A15.42}
\times~_2F_1\left[\begin{array}{ll}
2,1-s;\\
& 1-\frac{1}{\lambda^2}\\
~~2-s;\end{array}\right]-\frac{\Gamma(1+s)\Gamma(1-s)}{2},
\end{equation}
where $|1-\frac{1}{\lambda^2}|<1$ and $\mathfrak{R}(s)<1$.\\

\noindent
After further simplification, we arrive at the result (\ref{A15.22}).\\

\noindent
Throughout the paper, the semi-axes of the hyperboloid of one sheet $\frac{x^2}{a^2}+\frac{y^2}{b^2}-\frac{z^2}{c^2}=1$ are assumed in the form of $a>b>0$ and $c>0$.

\section{\bf  Closed form for the curved surface area of a hyperboloid of one sheet (Juggler's damroo of morha chair)}
{\bf Theorem 4.} The curved surface area of a hyperboloid of one sheet (whose axis is positive direction of $z$-axis) i.e, $z=c\left( \frac{x^2}{a^2}+\frac{y^2}{b^2}-1\right)^{\frac{1}{2}};~a>b>0\text{ and }
c>0$ bounded by the two planes $z=0$ ($x$-$y$ plane) and $z=H>0$ (a plane $||$ to $x$-$y$ plane) is given by
\begin{equation*}
\hat{S}~=\frac{ 2b^2c\sqrt{(\lambda^2-1)}~\pi}{(\lambda a)}\times~~~~~~~~~~~~~~~~~~~~~~~~~~~~~~~~~~~~~~~~~~~~~~~~~~~~~~~~~~~~~`
\end{equation*}
\begin{equation}\label{A15.100}
\times F_{2\text{{\bf:}}1;0;0}^{1\text{{\bf:}}1;1;2}\left( \begin{array}{ll}
~~~~~~~~~~~~~~[\frac{1}{2}:0,1,1]~~~~\text{\bf :}[\frac{1}{2}:1];[{2}:1];[\frac{-1}{2}:1],[\frac{-1}{2}:1]~\text{\bf;} \\
& \frac{b^2}{a^2}-1\text{\bf,}1-\frac{1}{\lambda^2}\text{\bf,}\frac{a^2(1-\lambda^2)}{c^2\lambda^2} \\

[\frac{-1}{2}:-1,0,1],[\frac{3}{2}:0,1,1]\text{\bf :}[1:1];~~-~~;~~~\text{------------}~~~~~~\text{\bf;}\end{array}\right)
\end{equation}
where $\lambda^2=\left(1+ \frac{H^2}{c^2}\right);~\lambda>1;~|1-\frac{1}{\lambda^2}|<1;~|\frac{b^2}{a^2}-1|<1$ and $|\frac{a^2(1-\lambda^2)}{c^2\lambda^2}|<1$.\\

\noindent
{\em{Remark:}} The above formula (\ref{A15.100}) is verified numerically through {\em Mathematica program}.\\
\begin{figure}[ht!]

\noindent\begin{minipage}{0.70\textwidth}% adapt widths of minipages to your needs
\includegraphics[width=\linewidth]{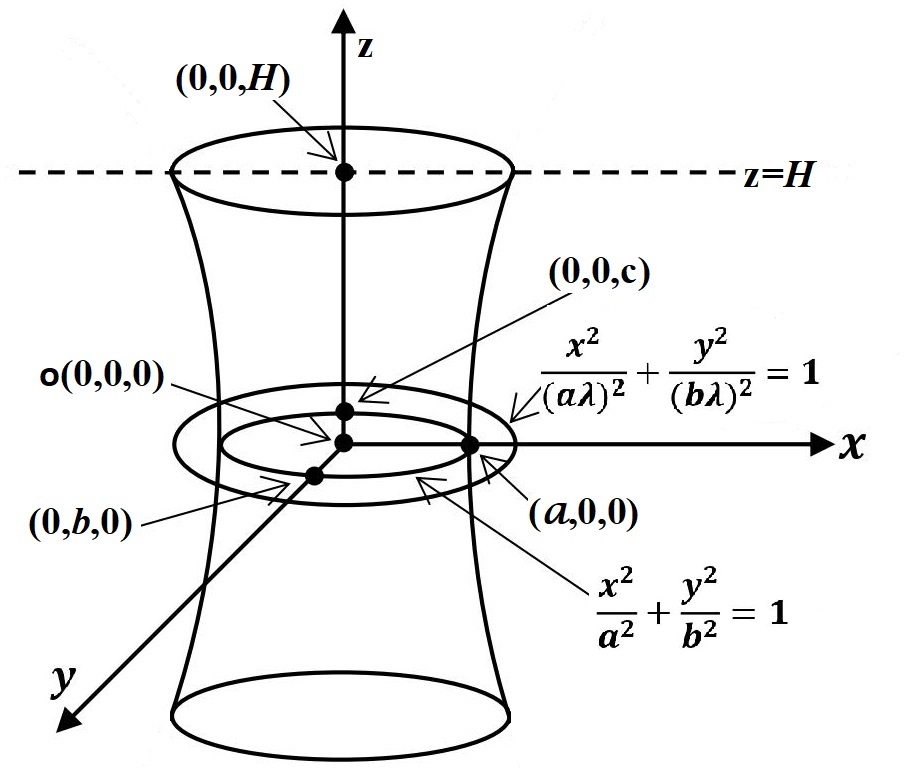}
\end{minipage}%
\hfill%
\caption[short]{Hyperboloid of one sheet $\frac{x^2}{a^2}+\frac{y^2}{b^2}-\frac{z^2}{c^2}=1$}
\end{figure}

\begin{proof}
\noindent
Equation of a hyperboloid of one sheet
\begin{equation}\label{A15.101}
\frac{x^2}{a^2}+\frac{y^2}{b^2}-\frac{z^2}{c^2}=1,~~a>b>0;c>0.~~~~~~~~~~~~~~~~~~~~
\end{equation}
Therefore,
\begin{equation}\label{A15.102}
\frac{x^2}{a^2}+\frac{y^2}{b^2}=\frac{z^2}{c^2}+1.~~~~~~~~~~~~~~~~~~~~~~~~~~~~~~~~~~~~
\end{equation}
The intersection of the surface (\ref{A15.102}) with the plane $z=H$, where $0<c< H $ (parallel to $x$-$y$ plane, lying above $x$-$y$ plane) will be an ellipse (in space) lying in the plane $z=H$ as well as lying on the surface of a hyperboloid of one sheet. The projection of that ellipse in $x$-$y$ plane will be

\begin{equation}\label{A15.103}
\frac{x^2}{a^2}+\frac{y^2}{b^2}=\frac{H^2}{c^2}+1=\lambda^2,~~~~~~~~~~~~~~~~~~~~~~~~~~~~~
\end{equation}
where $\lambda=\sqrt{\left( \frac{H^2}{c^2}+1\right) }$.\\

Again from equation (\ref{A15.102}) the equation of upper portion of hyperboloid of one sheet will be
\begin{equation}\label{A15.104}
z=c\left( \frac{x^2}{a^2}+\frac{y^2}{b^2}-1\right)^{\frac{1}{2}} ~~~~~~~~~~~~~~~~~~~~~~~~~~~~~~~~~~~~~~~~~~~
\end{equation}
\begin{equation}
\frac{\partial z}{\partial x}=\frac{cx}{a^2\sqrt{\left( \frac{x^2}{a^2}+\frac{y^2}{b^2}-1\right)}} ~~~~~~~~~~~~~~~~~~~~~~~~~~~~~~~~~~~~~~~~~~~
\end{equation}
\begin{equation}
\frac{\partial z}{\partial y}=\frac{cy}{b^2\sqrt{\left( \frac{x^2}{a^2}+\frac{y^2}{b^2}-1\right)}} ~~~~~~~~~~~~~~~~~~~~~~~~~~~~~~~~~~~~~~~~~~~
\end{equation}

\noindent
The projection of curved surface area of a hyperboloid of one sheet (\ref{A15.104}) in $x$-$y$ plane will be the area between the two concentric ellipses, given by\\
\begin{figure}[ht!]
\noindent\begin{minipage}{0.50\textwidth}% adapt widths of minipages to your needs
\includegraphics[width=\linewidth]{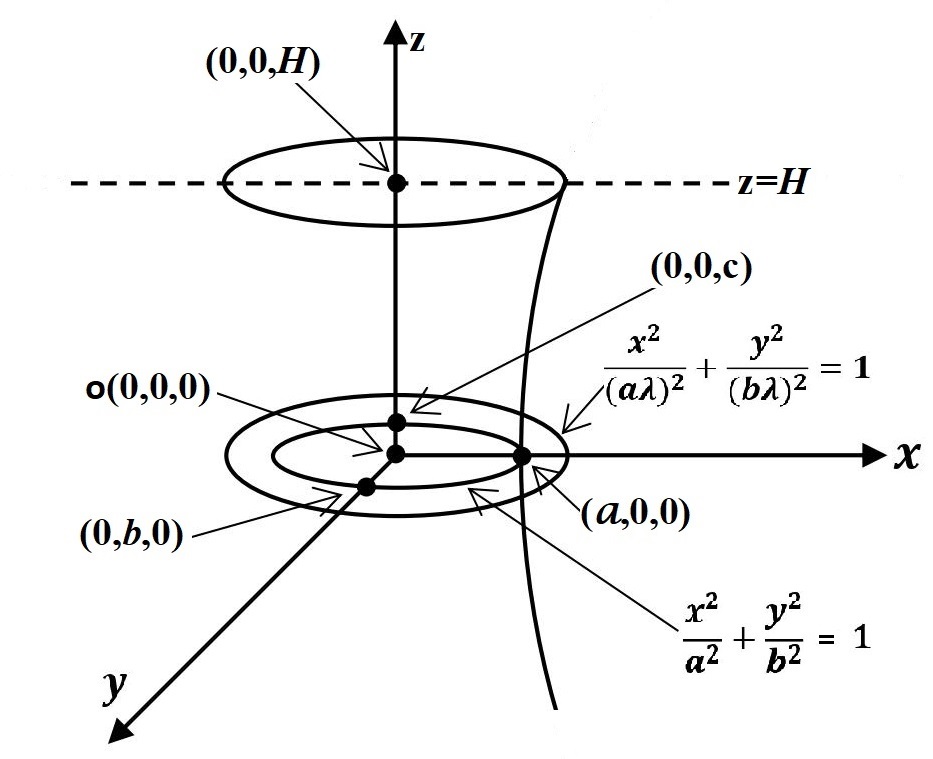}
\end{minipage}%
\hfill%
\caption[short]{Concentric ellipses}
\end{figure}

\begin{equation}\label{A15.43}
\frac{x^2}{(a\lambda)^2}+\frac{y^2}{(b\lambda)^2}=1,~~~~~~~~~~~~~~~~~~~~~~~~~~~~~
\end{equation}
\begin{equation}\label{A15.44}
\frac{x^2}{a^2}+\frac{y^2}{b^2}=1,~~~~~~~~~~~~~~~~~~~~~~~~~~~~~
\end{equation}
where $a>b $.\\

\noindent
Substitute the values of $\frac{\partial z}{\partial x}$ and $\frac{\partial z}{\partial y}$ in equation (\ref{a12.9}).
Therefore curved surface area of a hyperboloid of one sheet will be
\begin{equation}
\hat{S}~=\underbrace{\int\int}_{\stackrel{\text{over the area between two concentric}}{\text{   ellipses (\ref{A15.43}) and (\ref{A15.44})}}}\sqrt{\left\lbrace 1+\frac{c^2x^2}{a^4\left( \frac{x^2}{a^2}+\frac{y^2}{b^2}-1\right)}+\frac{c^2y^2}{b^4\left( \frac{x^2}{a^2}+\frac{y^2}{b^2}-1\right)}\right\rbrace }~dx ~dy.
\end{equation}
\begin{figure}[ht!]
\noindent\begin{minipage}{0.50\textwidth}% adapt widths of minipages to your needs
\includegraphics[width=\linewidth]{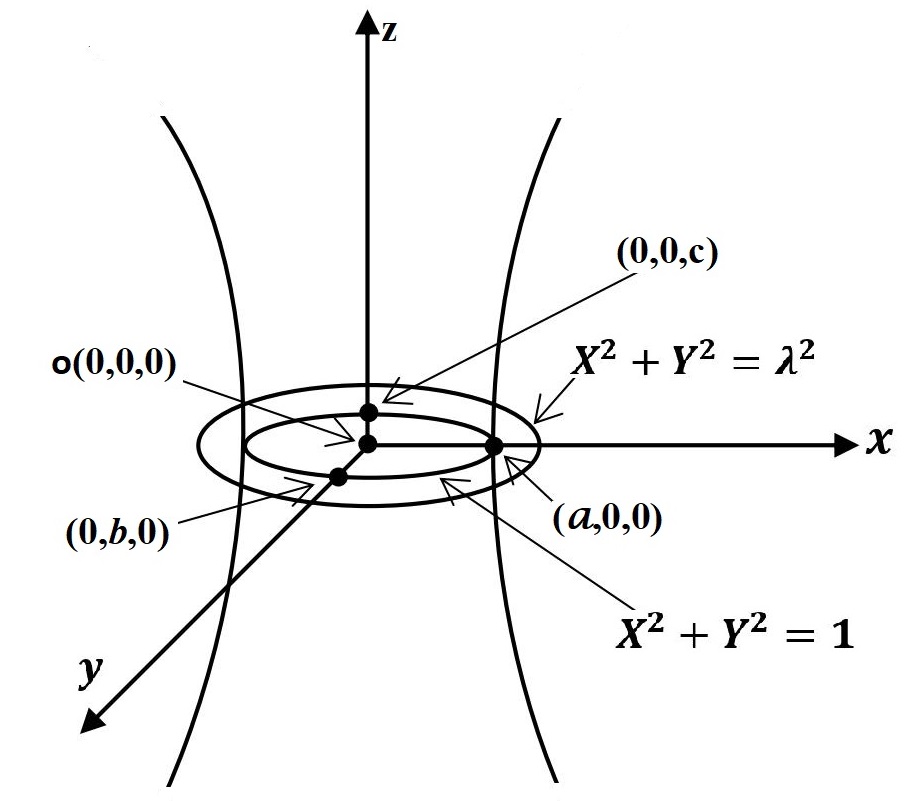}
\end{minipage}%
\hfill%
\caption[short]{Concentric circles}
\end{figure}
Put $x={a}{X},~y={b}{Y}$, then
\begin{equation}
\hat{S}=ab~\underbrace{\int\int}_{\stackrel{\text{over the area between two concentric}}{\text{  circles}~X^2+Y^2=\lambda^2,X^2+Y^2=1}}\sqrt{\left\lbrace 1+\frac{c^2}{\left( X^2+Y^2-1\right)}\left( \frac{X^2}{a^2}+\frac{Y^2}{b^2}\right) \right\rbrace} ~dX ~dY
\end{equation}

Put $X=r\cos \theta,~Y=r\sin \theta$, then $dX~dY=r~dr~d\theta$.
\begin{equation}\label{A15.2000}
\text{Therefore }~\hat{S}=ab~\int_{\theta=-\pi}^{\pi}\left( \int_{r=1}^{\lambda}\left\lbrace 1+\frac{c^2r^2}{\left( r^2-1\right)}\left( \frac{\cos^2\theta}{a^2}+\frac{\sin^2\theta}{b^2}\right) \right\rbrace^{\frac{1}{2}} ~rdr\right)  ~d\theta
\end{equation}
{\bf Remark: }Since we have no standard formula of definite/indefinite integrals in the literature of integral calculus for the integration with respect to $"r"$ and $"\theta"$ in double integral (\ref{A15.2000}). Therefore we can solve such integrals exactly through hypergeometric function approach.
\begin{equation}\label{A15.45}
\text{Therefore }~\hat{S}=ab~\int_{\theta=-\pi}^{\pi}\int_{r=1}^{\lambda}{}_{1}F_{0}\left[\begin{array}{ll}
\frac{-1}{2};\\
& \frac{c^2r^2}{\left( 1-r^2\right)}\left( \frac{\cos^2\theta}{a^2}+\frac{\sin^2\theta}{b^2}\right) \\
~-;\end{array}\right]~rdr ~d\theta.
\end{equation}
Since there is uncertainty about the argument of$~_1F_0$ in equation (\ref{A15.45}), because the argument of $~_1F_0$ in equation (\ref{A15.45}) may be greater than 1. Therefore applying contour integral (\ref{D12.100}) of ${}_{1}F_{0}(.)$ in equation (\ref{A15.45}), we get
\begin{equation}
\hat{S}=ab~\int_{\theta=-\pi}^{\pi}\int_{r=1}^{\lambda}\left[ \frac{1}{2\pi i \Gamma\left(\frac{-1}{2}\right) }~{\int}_{s=-i\infty}^{+i\infty}~\Gamma(-s)\Gamma( \frac{-1}{2}+s) \left\lbrace \frac{-c^2r^2}{(1-r^2)}\left( \frac{\cos^2\theta}{a^2}+\frac{\sin^2\theta}{b^2}\right) \right\rbrace^s ds \right]  ~rdr~d\theta,
\end{equation}
where $|\arg\left\lbrace \frac{c^2r^2}{(r^2-1)}\left( \frac{\cos^2\theta}{a^2}+\frac{\sin^2\theta}{b^2}\right) \right\rbrace|<\pi$.
\begin{equation*}
\text{Therefore }\hat{S}=\frac{-ab}{4\pi \sqrt{\pi}~i}~{\int}_{s=-i\infty}^{+i\infty}~ \Gamma(-s)\Gamma( \frac{-1}{2}+s) (-c^2)^s\left\lbrace ~\int_{\theta=-\pi}^{\pi}\left( \frac{\cos^2\theta}{a^2}+\frac{\sin^2\theta}{b^2}\right) ^s~d\theta\right\rbrace\times
\end{equation*}

\begin{equation}\label{A15.46}
\times \left\lbrace \int_{r=1}^{\lambda}\frac{r^{2s+1}}{(1-r^2)^s}~dr \right\rbrace  ds.
\end{equation}
Since $a>b$ and $\lambda>1$, therefore employing the formulae (\ref{G12.10}) and (\ref{A15.22}) in (\ref{A15.46}), we get
\begin{equation*}
\hat{S}= \frac{-b^2(1-\lambda^2)}{4\lambda^2 \sqrt{\pi}~i}~{\int}_{s=-i\infty}^{+i\infty}\frac{
 \Gamma(-s)\Gamma( \frac{-1}{2}+s)\Gamma(s-1)}{\Gamma(s)(1-\lambda^2)^s}\left( \frac{-c^2\lambda^2}{a^2}\right)^s~{}_{2}F_{1}\left[\begin{array}{ll}
\frac{1}{2},1+s;\\
& 1-\frac{b^2}{a^2}\\
~1;\end{array}\right]\times
\end{equation*}

\begin{equation}
\times~{}_{2}F_{1}\left[\begin{array}{ll}
2,1-s;\\
& 1-\frac{1}{\lambda^2}\\
~~2-s;\end{array}\right]~ds,
\end{equation}
where $(1-\frac{b^2}{a^2})<1$ and $(1-\frac{1}{\lambda^2})<1$.
\begin{equation*}
\hat{S}=  \frac{-b^2(1-\lambda^2)}{4\lambda^2 \sqrt{\pi}~i}{\int}_{s=-i\infty}^{+i\infty}\frac{
 \Gamma(-s)\Gamma( \frac{-1}{2}+s)\Gamma(s-1)}{\Gamma(s)(1-\lambda^2)^s}\left( \frac{-c^2\lambda^2}{a^2}\right)^s\sum_{m=0}^{\infty}\frac{\left( \frac{1}{2}\right)_m(1+s)_m\left( 1-\frac{b^2}{a^2}\right)^m  }{(1)_m~m!}\times
 \end{equation*}
 \begin{equation}
 \times\sum_{n=0}^{\infty}\frac{(2)_n(1-s)_n~\left( 1-\frac{1}{\lambda^2}\right) ^n}{(2-s)_n~n!}~ds.
\end{equation}
Interchanging the order of summation and integration, we get

\begin{equation*}
\hat{S}= \frac{b^2(1-\lambda^2)}{4\lambda^2 \sqrt{\pi}~i}~\sum_{m,n=0}^{\infty}\frac{\left(\frac{1}{2}\right)_m(2)_n\left( 1-\frac{b^2}{a^2}\right) ^m\left( 1-\frac{1}{\lambda^2}\right) ^n }{(1)_m~m!~n!}\times~~~~~~~~~~~~~~~~~~~~~~~~~~~~
\end{equation*}
\begin{equation}
\times{\int}_{s=-i\infty}^{+i\infty}~
 \frac{\Gamma(-s)\Gamma( \frac{-1}{2}+s)\Gamma( 1+m+s)\Gamma( 1+n-s)}{\Gamma(1+s)\Gamma( 2+n-s)}\left( \frac{-c^2\lambda^2}{a^2(1-\lambda^2)}\right)^s~ds.
\end{equation}
Applying the definition (\ref{E12.12}) of $G$-function, we get

\begin{equation}\label{A15.47}
\hat{S}=\frac{b^2(1-\lambda^2)\sqrt{\pi}}{2\lambda^2 }~\sum_{m,n=0}^{\infty}\frac{\left(\frac{1}{2}\right)_m(2)_n\left( 1-\frac{b^2}{a^2}\right) ^m\left( 1-\frac{1}{\lambda^2}\right) ^n }{(1)_m~m!~n!}~G^{2~2}_{3~3}\left(\frac{c^2\lambda^2}{a^2(\lambda^2-1)}\left|   \begin{array}{ll}
~\frac{3}{2},-m~;2+n\\
0,1+n;~0  \end{array}\right.\right),
\end{equation}
where $|\frac{c^2\lambda^2}{a^2(\lambda^2-1)}|>1$.\\

Employing the decomposition formula for Meijer's $G$- function (\ref{A15.21}) in equation (\ref{A15.47}), we get

\begin{equation*}
\hat{S}= \frac{2b^2c\sqrt{(\lambda^2-1)}~{\pi}}{\lambda a}\sum_{m,n=0}^{\infty}\frac{\left(\frac{1}{2}\right)_m\left(\frac{3}{2}\right)_m(2)_n\left(\frac{1}{2}\right)_n\left( 1-\frac{b^2}{a^2}\right) ^m\left( 1-\frac{1}{\lambda^2}\right) ^n }{(1)_m\left(\frac{3}{2}\right)_n~m!~n!}\times
\end{equation*}
\begin{equation}
\times~_{3}F_{2}\left[\begin{array}{ll}
~~\frac{-1}{2},\frac{-1}{2},\frac{1}{2}+n;\\
& \frac{a^2(1-\lambda^2)}{c^2\lambda^2}\\
\frac{-1}{2}-m,\frac{3}{2}+n;\end{array}\right],
\end{equation}
where $|\frac{a^2(1-\lambda^2)}{c^2\lambda^2}|<1$.
\begin{equation*}
\text{Therefore }\hat{S}= \frac{2b^2c\sqrt{(\lambda^2-1)}~{\pi}}{\lambda a}\sum_{m,n=0}^{\infty}\frac{\left(\frac{1}{2}\right)_m\left(\frac{3}{2}\right)_m(2)_n\left(\frac{1}{2}\right)_n\left( 1-\frac{b^2}{a^2}\right) ^m\left( 1-\frac{1}{\lambda^2}\right) ^n }{(1)_m\left(\frac{3}{2}\right)_n~m!~n!}\times
\end{equation*}

\begin{equation}\label{A15.48}
\times~\sum_{p=0}^{\infty}\frac{\left( \frac{-1}{2}\right)_p\left( \frac{-1}{2}\right)_p\left( \frac{1}{2}+n\right)_p\left( \frac{a^2(1-\lambda^2)}{c^2\lambda^2}\right)^p}{
\left( \frac{-1}{2}-m\right)_p \left( \frac{3}{2}+n\right)_p~p!}
\end{equation}
\begin{equation}\label{A15.49}
= \frac{2b^2c\sqrt{(\lambda^2-1)}~{\pi}}{\lambda a}\sum_{m,n,p=0}^{\infty}\frac{\left(\frac{1}{2}\right)_{n+p}\left(\frac{1}{2}\right)_m(2)_n\left(\frac{-1}{2}\right)_p\left(\frac{-1}{2}\right)_p\left( \frac{b^2}{a^2}-1\right) ^m\left( 1-\frac{1}{\lambda^2}\right) ^n\left( \frac{a^2(1-\lambda^2)}{c^2\lambda^2}\right)^p }{\left( \frac{-1}{2}\right)_{-m+p}\left( \frac{3}{2}\right)_{n+p}(1)_m~m!~n!~p!}.
\end{equation}
Expressing the result (\ref{A15.49}) in terms of Srivastava and Daoust hypergeometric function (\ref{AB15.42}), we arrive at the result (\ref{A15.100}).
\end{proof}

\section{\bf Volume of a hyperboloid of one sheet}
The volume of a hyperboloid of one sheet (Juggler's damroo or morha chair) having the vertical height $H$, $a$ and $b$, the lengths of semi-major axis and semi-minor axis respectively, of the smaller elliptic base, is given by

\begin{equation}\label{AB15.125}
V=\pi ab H \left( 1+\frac{H^2}{3c^2}\right),  ~~~~~~~~~~~~~~~~~~~~~~~~~~
\end{equation}
\begin{equation}
\text{where}~ c^2=\frac{H^2}{(\lambda^2-1)}~~~~~~~~~~~~~~~~~~~~~~~~~~~~~~~~~~~~~~~~~~~~~~~
\end{equation}
\begin{equation}\label{ABC15.125}
\text{and }\lambda=\frac{\text{Length of semi-major (or minor) axis of larger elliptic base}}{\text{Length of semi-major (or minor) axis of smaller elliptic base}}.
\end{equation}
The equation (\ref{ABC15.125}) is one of the important condition for hyperboloid of one sheet. In this connection see figure 5.\\

\noindent
Note: The lengths of semi-major and semi-minor axes of larger and smaller elliptic bases of given model (Damroo or morha chair) can be easily measured and after that $\lambda$ and $c$ are calculated easily.

\begin{figure}[ht!]

\noindent\begin{minipage}{0.60\textwidth}% adapt widths of minipages to your needs
\includegraphics[width=\linewidth]{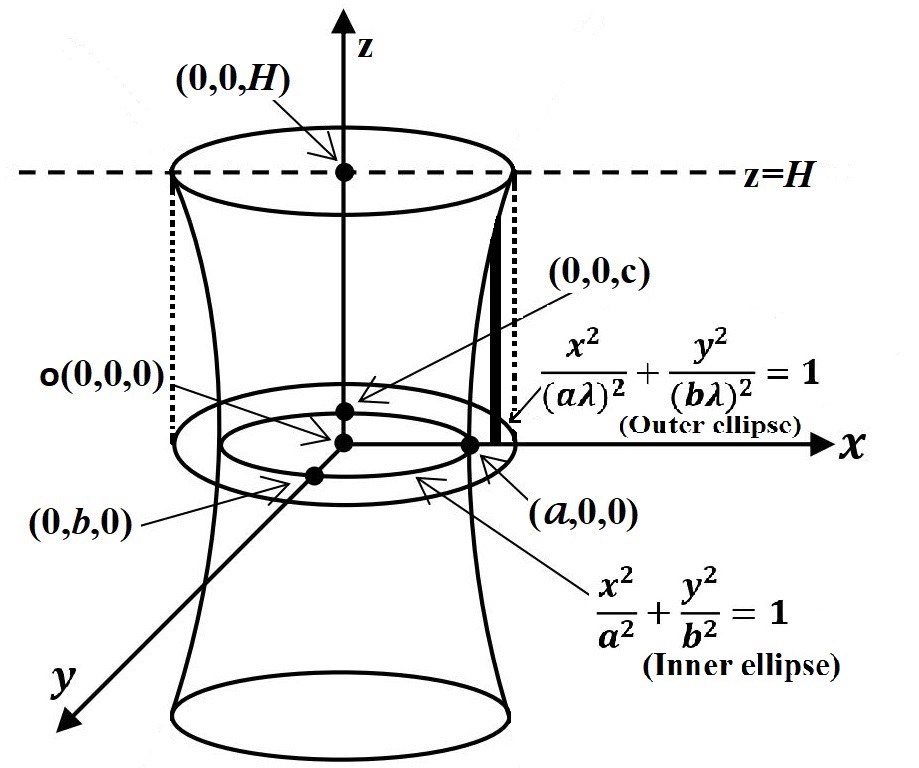}
\end{minipage}%
\hfill%
\caption[short]{Hyperboloid of one sheet $\frac{x^2}{a^2}+\frac{y^2}{b^2}-\frac{z^2}{c^2}=1$}
\end{figure}

\begin{proof} 
\noindent
From the equation of a hyperboloid of one sheet, we have
\begin{equation}
z=c\left( \frac{x^2}{a^2}+\frac{y^2}{b^2}-1\right)^{\frac{1}{2}}; ~~a>b>0\text{ and}~c>0\text{ taking along positive } z \text{ axis.}
\end{equation}
Now volume of a hyperboloid of one sheet bounded by the two planes $z=0$ and $z=H$ is $V=V_c-V_b$, where $V_c$ is the volume of a right elliptic cylinder whose base is an ellipse having the semi-axes $a\lambda$ and $b\lambda$ and vertical height $H$ i.e, ($V_c=\pi ab\lambda^2H$) and $V_b$ is the volume lying {\bf between} outside of the surface of a hyperboloid of one sheet and inside the curved surface of right elliptic cylinder.
\begin{equation}
\text{Then }V_b~= \underbrace{\int\int}_{\stackrel{\text{over the projection of upper portion of hyperboloid}}{\text{ of one sheet bounded by $z=0$ and $z=H$} }}\int_{z=0}^{z=\alpha}dz~dy ~dx;~~\alpha=c\sqrt{\left( \frac{x^2}{a^2}+\frac{y^2}{b^2}-1\right)}~~~~~~~~~~~~~~~~
\end{equation}
and $\lambda =\sqrt{\frac{H^2}{c^2}+1}.$
\begin{equation}
\text{Therefore }V_b=\underbrace{\int\int}_{\stackrel{\text{over the area between the ellipses}}{ \frac{x^2}{(a\lambda)^2}+\frac{y^2}{(b\lambda)^2}=1\text{ and }\frac{x^2}{a^2}+\frac{y^2}{b^2}=1}} c~\sqrt{\left( \frac{x^2}{a^2}+\frac{y^2}{b^2}-1\right)}~dy ~dx~~~~~~~~~~~~~~~~~~~
\end{equation}
Put $x=aX,~y=bY$, then
\begin{equation}
{V_b}~=abc\underbrace{\int\int}_{\stackrel{\text{over the area between the circles}}{X^2+Y^2=\lambda^2\text{ and }X^2+Y^2=1} } \sqrt{\left( X^2+Y^2-1\right)}~dY ~dX.~~~~~~~~~~~~~~~~~~~
\end{equation}

Take $X=r\cos \theta,~Y=r\sin \theta$, then
\begin{equation*}
{V_b}~=abc~\left( \int_{\theta=-\pi}^{\pi}~d\theta\right) \left(  \int_{r=1}^{\lambda}r\sqrt{\left( r^2-1\right)}dr \right) 
\end{equation*}
\begin{equation*}
=\frac{2 \pi abc}{3}\left[ (r^2-1)^{\frac{3}{2}}\right]_1^\lambda
\end{equation*}

\begin{equation}\label{A15.144}
V_b=\frac{2\pi abc(\lambda^2-1)^{\frac{3}{2}}}{3}.
\end{equation}
On taking the difference of $V_c$ and $V_b$, we get
\begin{equation}\label{A15.125}
V=\pi ab \lambda^2H-\frac{2\pi abc\left(\lambda^2-1\right)^{\frac{3}{2}}}{3} ;~~0<c<H,
\end{equation}
where $\lambda\left( =\sqrt{\frac{H^2}{c^2}+1}\right)>1$.\\
Substituting the value of $\lambda$ in equation (\ref{A15.125}) and simplifying further, we arrive at the required volume (\ref{AB15.125}) of a hyperboloid of one sheet $z=c\left( \frac{x^2}{a^2}+\frac{y^2}{b^2}-1\right)^{\frac{1}{2}};a>b>0$ and $c>0 $ bounded by the two planes $z=0$ and $z=H>0$.

\end{proof}

\section{\bf Conclusion}
In this paper, we obtained the closed form for the exact curved surface area of a hyperboloid of one sheet $\frac{x^2}{a^2}+\frac{y^2}{b^2}-\frac{z^2}{c^2}=1$ intercepted by the two parallel planes $z=0$ and $z=H$, through hypergeometric function approach i.e, by using series rearrangement technique, Mellin-Barnes type contour integral representations of generalized hypergeometric function$~_pF_q(z)$, Meijer's $G$-function and decomposition formula for Meijer's $G$-function; in terms of generalized hypergeometric function of Srivastava and Daoust. The formula is neither available in the literature of mathematics nor found in any mathematical tables. Moreover, we also derived the formula for the volume of a hyperboloid of one sheet.
We conclude that many formulas for curved surface areas of other three dimensional figures can be derived in an analogous manner, using Mellin-Barnes contour integration. Moreover, the results deduced above (presumably new), have potential applications in the fields of applied mathematics, statistics and engineering sciences.\\

\noindent
{\bf Conflicts of interests:} The authors declare that there are no conflicts of interest.


\begin{thebibliography}{999}
\normalsize
\bibitem{Abramowitz} Abramowitz, M. and Stegun, I.A.; {\em Handbook of Mathematical Functions with Formulas, Graphs and Mathematical Tables}, Reprint of the 1972 Edition, Dover Publications, Inc., New York, 1992.
\bibitem{And1999} Andrews, G.E., Askey, R. and Roy, R.; {\em Special Functions}, Cambridge University Press, Cambridge, UK, 1999.

\bibitem{Erd1} Erd\'{e}lyi, A., Magnus,  W., Oberhettinger,  F.  and  Tricomi, F.G.; {\em Higher Transcendental Functions}, Vol. I, McGraw-Hill Book Company, New York, Toronto and London, 1955.

\bibitem{Lauricella93158} Lauricella, G.; Sulle funzioni ipergeometriche a più variabili, {\em Rend. Circ. Mat. Palermo}, {\bf 7} (1893), 111-158.

\bibitem{Lebedev1972} Lebedev, N.N.; {\em Special Functions and their Applications,} translated by Richard A. Silverman. Prentice-hall, Inc, Englewood Cliffs, N.J., 1965.

\bibitem{Luke169} Luke, Y.L.; {\em The Special Functions and Their Approximations}, Vol. I, Academic Press, 1969.

\bibitem{MagOberSoni1966} Magnus, W., Oberhettinger, F. and Soni, R.P.; {\em Some Formulas and Theorems for the Special Functions of Mathematical Physics}, Third Enlarged Edition, Springer-Verlag, New York, 1966.

 \bibitem{Mathaisaxena} Mathai, A.M. and Saxena, R.K.; {\em Lecture notes in Mathematics No.348: Generalized hypergeometric functions with Applications in statistics and physical sciences}, Springer-Verlag, Berlin Heidelberg, New York, 1973.
 
\bibitem{Prudbrimari90}	Prudnikov, A.P., Brychknov, Yu. A. and Marichev, O.I.; {\em Integrals and Series}, Vol. III: {More special functions}, Nauka Moscow, 1986 (in Russian);(Translated from the Russian by G.G.Gould), Gordon and Breach Science Publishers, New York, Philadelphia London, Paris, Montreux, Tokyo, Melbourne, 1990.

\bibitem{Rain} Rainville, E.D.;  {\em Special Functions}, The Macmillan Co. Inc., New York 1960; Reprinted by Chelsea publ. Co., Bronx, New York, 1971.

\bibitem{Slater1966} Slater, L.J.; {\em Generalized Hypergeometric Functions}, Cambridge Univ., Press, New York, 1966.

\bibitem{Sridaoust1969202} Srivastava, H.M. and Daoust, M.C.; On Eulerian integrals associated with Kamp${\acute e}$ de F${\acute e}$riet's function, \textit{Publ. Inst. Math. (Beograd) (N.S.)}, {\bf 9} (23) (1969), 199-202.
\bibitem{Sridaoust1969457}Srivastava, H.M. and Daoust, M.C.; Certain generalized Neumann expansions associated with the Kamp\'{e} de F\'{e}riet's function, \textit{Nederl. Akad. Wetensch. Proc. Ser. A, 72=Indag. Math.,} {\bf 31} (1969), 449-457.
\bibitem{Sridaoust1972} Srivastava, H.M. and Daoust, M.C.; A note on the convergence of Kamp\'{e} de F\'{e}riet's Double hypergeometric series, \textit{Math. Nachr.,} {\bf 53} (1972), 151-159.

\bibitem{SriMano} Srivastava, H.M. and  Manocha, H.L.; {\em A Treatise on Generating Functions}, Halsted Press (Ellis Horwood Limited, Chichester), John Wiley and Sons, New York, Chichester, Brisbane and Toronto, 1984.

\bibitem{Wikipedia} Wolfram research, Meijer $G$- function: specific values (subsection 03/01), {\em http://functions.wolfram.com/HypergeometricFunctions/MeijerG/03/01/Show\\All.html}, 1-99.

\bibitem{Wright1935} Wright, E.M.; The asymptotic expansion of the generalized hypergeometric function, {\em J. London  Math. Soc.}, {\bf 10} (1935), 286-293.
\bibitem{Wright1940} Wright, E.M.; The asymptotic expansion of the generalized hypergeometric function,  {\em Proc. London  Math. Soc.},{\bf 46} (2) (1940), 389-408.
\end{thebibliography}
\end{document}